\documentclass[12pt]{article}
\usepackage{amsmath}
\usepackage{amssymb}
\usepackage{hyperref}
\usepackage{fullpage}
\usepackage{booktabs}
\usepackage{xspace}
\usepackage{longtable}
\usepackage{authblk}
\usepackage{xcolor}
\usepackage{graphicx}
\usepackage{microtype}
\usepackage[toc,page]{appendix}

\hypersetup{
    colorlinks=true,
    linkcolor={red!50!black},
    citecolor={blue!50!black},
    urlcolor={blue!80!black}
}

\newcommand{\OEISeditsMain}{141\xspace}
\newcommand{\OEISeditsExtended}{6\xspace}

\newcommand{\OEIS}[1]
{\href{https://oeis.org/#1}{\texttt{#1}}}

\newcommand{\SEQ}{S}
\newcommand{\QUERY}{\mathcal{Q}}
\newcommand{\CLASS}{\mathcal{C}}
\newcommand{\ie}[0]{i.e.\ }
\newcommand{\eg}[0]{e.g.\ }

\newcommand{\VARautomorphismgroupn}{| \text{Aut}(g) |}
\newcommand{\VARedgeconnectivity}{\kappa_{E}(g)}
\newcommand{\VARvertexconnectivity}{\kappa_{V}(g)}
\newcommand{\VARdiameter}{\text{dia}(g)}
\newcommand{\VARgirth}{\text{girth}(g)}
\newcommand{\VARnarticulationpoints}{a(g)}
\newcommand{\VARmaximalindependentvertexset}{\alpha_V(g)}
\newcommand{\VARmaximalindependentedgeset}{\alpha_E(g)}

\newcommand{\VARchromatic}{\chi}
\newcommand{\VARchromaticpolynomial}[1]{P_g(#1)}
\newcommand{\VARfractionalchromaticnumber}[1]{\VARchromatic_{#1}(g)}
\newcommand{\VARchromaticnumber}{\VARchromatic(g)}

\newcommand{\VARsubgraph}{s}
\newcommand{\namedsubgraph}[1]{\VARsubgraph{}(g,#1)}
\newcommand{\VARissubgraphfreeKthree}{\namedsubgraph{K_3}}
\newcommand{\VARissubgraphfreeKfour}{\namedsubgraph{K_4}}
\newcommand{\VARissubgraphfreeKfive}{\namedsubgraph{K_5}}
\newcommand{\VARissubgraphfreeCfour}{\namedsubgraph{C_4}}
\newcommand{\VARissubgraphfreeCfive}{\namedsubgraph{C_5}}

\newcommand{\subgraphBULL}{B_5}
\newcommand{\subgraphDIAMOND}{D_4}
\newcommand{\subgraphBOWTIE}{\bowtie}
\newcommand{\subgraphOPENBOWTIE}{\rtimes}
\newcommand{\VARissubgraphfreebull}{\namedsubgraph{\subgraphBULL}}
\newcommand{\VARissubgraphfreediamond}{\namedsubgraph{\subgraphDIAMOND}}
\newcommand{\VARissubgraphfreeopenbowtie}{\namedsubgraph{\subgraphBOWTIE}}
\newcommand{\VARissubgraphfreebowtie}{\namedsubgraph{\subgraphOPENBOWTIE}}

\newcommand{\VARfractionaldualitygapvertexchromatic}{X}
\newcommand{\VARdistanceregular}{D}
\newcommand{\VARhamiltonian}{H}
\newcommand{\VARbipartite}{B}
\newcommand{\VAReulerian}{E}
\newcommand{\VARplanar}{P}
\newcommand{\VARtree}{T}
\newcommand{\VARchordal}{C}
\newcommand{\VARkregular}{R}
\newcommand{\VARstronglyregular}{R^*}
\newcommand{\VARintegral}{I}
\newcommand{\VARrealspectrum}{I^*}

\newcommand{\indicatorfunctionX}[1]{{#1}(g)}
\newcommand{\VARhasfractionaldualitygapvertexchromatic}
{\indicatorfunctionX{\VARfractionaldualitygapvertexchromatic}}
\newcommand{\VARisdistanceregular}
{\indicatorfunctionX{\VARdistanceregular}}
\newcommand{\VARishamiltonian}
{\indicatorfunctionX{\VARhamiltonian}}
\newcommand{\VARisbipartite}
{\indicatorfunctionX{\VARbipartite}}
\newcommand{\VARiseulerian}
{\indicatorfunctionX{\VAReulerian}}
\newcommand{\VARisplanar}
{\indicatorfunctionX{\VARplanar}}
\newcommand{\VARistree}
{\indicatorfunctionX{\VARtree}}
\newcommand{\VARischordal}
{\indicatorfunctionX{\VARchordal}}
\newcommand{\VARiskregular}
{\indicatorfunctionX{\VARkregular}}
\newcommand{\VARisstronglyregular}
{\indicatorfunctionX{\VARstronglyregular}}
\newcommand{\VARisintegral}
{\indicatorfunctionX{\VARintegral}}
\newcommand{\VARisrealspectrum}
{\indicatorfunctionX{\VARrealspectrum}}

\begin{document}

\title{Integer sequence discovery from small graphs}
\author[1]{Travis Hoppe}
\author[2]{Anna Petrone}
\affil[1]{National Institutes of Health, National Institute of Diabetes and Digestive and Kidney Diseases, Bethesda, Maryland}
\affil[2]{University of Maryland, Department of Civil Engineering}

\date{\today}
\maketitle

\begin{abstract}
We have exhaustively enumerated all simple, connected graphs of a finite order and have computed a selection of invariants over this set.
Integer sequences were constructed from these invariants and checked against the Online Encyclopedia of Integer Sequences (OEIS).
\OEISeditsMain new sequences were added and \OEISeditsExtended sequences were appended or corrected.
From the graph database, we were able to programmatically suggest relationships among the invariants.
It will be shown that we can readily visualize any sequence of graphs with a given criteria. 
The code has been released as an open-source framework for further analysis and the database was constructed to be extensible to invariants not considered in this work.
\end{abstract}

\section{Introduction}

There is a long history of public graph databases. 
Databases originally found only in print, such as the \textit{Atlas of Graphs}\cite{read1998atlas}, have rapidly expanded to the electronic medium.
These databases range from those of mathematical and algorithmic interest \cite{de2003large, brinkmann2013house}, to those cataloging structures found in the applied sciences such as ChemSpider\cite{pence2010chemspider}, RNA topologies\cite{gan2004rag, laing2013predicting} or social databases \cite{mcauley2014circles,brandes2009wikipedia,oestreicher2006amazon}.
Due to the rapid growth in the number of unique isomorphic graphs, the currently available databases are specialized in the number of graphs considered; 
a judicious choice often restricts the study to an interesting and more manageable subset.

We proceed with the assumption however, that \textit{a priori} all graphs could be interesting given the right question.
This is similar to the GraPHedron project\cite{melot2008facet}, which attempts to formulate conjectures by searching for graphs bounded by an inequality or constraint.
Our objective is more elementary. We aim to compute a large, comprehensive database of graphs and their respective invariants. 
Such a database will allow new forms of discovery, some of which will be directly explored in this paper.
For example, the integer sequences formed by the invariants can be systematically explored and compared to those already known.
These sequences, and the set of graphs that belong to them, can be used to explore a basic set of relations among the invariants. 
Since the input to GraPHedron consists of a set of invariants, a larger input set will amplify the predictive power.
Additionally, the creation of a large, centralized database will serve as a useful reference for benchmarking various algorithms.
Finally, a comprehensive database provides pedagogic value, as representative graphs from any considered sequence can be rapidly visualized.

To this end, we have created the \textit{Encyclopedia of Finite Graphs} (henceforth Encyclopedia), a database of invariants\cite{Travis2014SimpleConnectedGraphs} and the software to fully populate it\cite{Travis2014Encyclopdia}. 
The code and the database have been released under an open source license.
The intention is for new invariants to be an added to the project as various algorithms become available.

Once built, the Encyclopedia readily yields integer sequences formed by matching the number of graphs to an invariant constraint at each order.
Many such sequences have already been found and cataloged in another database, the Online Encyclopedia of Integer Sequences (OEIS)\cite{sloane2003line}.
The OEIS was created by Neil Sloane in 1964 as a graduate student during his studies of combinatorial problems.
Since then, the database has grown to over 250,000 sequences and is highly cited, with over 3,000 citations to date.
The sequences are of general interest, spanning topics such as number theory, combinatorics, and graph theory.
A given sequence may contain any number of terms, ranging from at least four up to as many as 500,000 (in the cases where the sequence admits a readily computable expression).
Collecting and storing integer sequences in one place allows a researcher, who perhaps comes across the first few values of an unknown sequence, to be able to quickly look up subsequent values. 
The OEIS not only provides the numerical values but seeks to function as a true encyclopedia, with cross references to related sequences, references to other literature, and formulas when known.
One of the primary goals of this paper is to systematically expand the sequences involving graph invariants known to the OEIS database.
Through our exhaustive enumeration of small graphs, we were able to submit \OEISeditsMain new sequences to the OEIS and extend \OEISeditsExtended existing sequences.

A graph invariant is any property that is preserved under isomorphism. 
Invariants can be simple binary properties (planarity), integers (automorphism group size), polynomials (chromatic polynomials), rationals (fractional chromatic numbers), complex numbers (adjacency spectra), sets (dominion sets) or even graphs themselves (subgraph and minor matching). 
We are primarily concerned with the sequences produced by graph invariants, \ie the combinatorial problem of how many graphs of a given class satisfy a particular criteria.
Let a graph be defined as the pair $G = (V,E)$, where $V$ is a set of vertices and $E$ is a set of edges. 
Define $\CLASS$ as a class which forms an isomorphically distinct set of graphs that satisfy a specified criteria.
Group the graphs into non-overlapping subsets such that
\begin{equation}
\CLASS = \CLASS_1 \cup \CLASS_2 \cup \CLASS_3 \cup \ldots
\end{equation}
where $\CLASS_n$ contains only graphs of order $n$.
From here, define an ordered sequence of subsets of the graph class
\begin{align}
\QUERY(f, \CLASS) 
&= f(\CLASS_1), f(\CLASS_2), f(\CLASS_3), \ldots
\end{align}
where $f$ is some invariant condition that selects from each set $\CLASS_n$.
Since $\QUERY$ selects from the graphs, we call $\QUERY$ the \textit{query}.

Let 
$\SEQ(f, \CLASS) = |f(\CLASS_1)|, |f(\CLASS_2)|, |f(\CLASS_3)|, \ldots$
be the sequence of integers defined by $\QUERY(f, \CLASS)$. 
For example, if $\VARistree$ is the $\{0,1\}$ indicator function that determines if the graph is a tree, and $\CLASS^\prime$ is the set of all simple unlabeled connected graphs, then
\begin{align}
S(\VARistree=1, \CLASS^\prime) = 1, 1, 1, 1, 2, 3, 6, 11, 23, \ldots
\end{align}
which is sequence \OEIS{A000055} in the OEIS.
This particular sequence is well-known and easily computable.
We have evaluated a range of invariants, from those that are computable in polynomial time, to some that are known to be $\mathcal{\#P}$-complete.
While a few sequences were merely extended, other invariants, such as the independence number, were hitherto unknown to OEIS and have produced novel sequences (\OEIS{A243781}-\OEIS{A243784}).

In addition to contributing to the OEIS, a secondary goal is to identify relationships between graph invariants.
Two queries of the same class are subsets of each other $\QUERY_a(f,\CLASS) \subseteq \QUERY_b(g, \CLASS)$, if $f(\CLASS_i) \subseteq g(\CLASS_i)$ for all $i\ge0$. 
Equality of two queries $\QUERY_a = \QUERY_b$, implies $\QUERY_a \subseteq \QUERY_b$ and $\QUERY_b \subseteq \QUERY_a$. 
We say that a relation between two invariants conditions is \textit{suggestive} to order $n$, $\QUERY_a \subseteq_n \QUERY_b$, if $f(\CLASS_i) \subseteq g(\CLASS_i)$ for $0 \le i \le n$.
They are \textit{exclusive} to order $n$, $\QUERY_a \cap_n \QUERY_b$, if $f(\CLASS_i) \cap g(\CLASS_i) = \emptyset$ for $0 \le i \le n$.
We can quickly filter candidate relations by noting that it is necessary but not sufficient for the same conditions to hold for sequences, \eg $\SEQ_a \ne \SEQ_b \implies  \QUERY_a \ne \QUERY_b$.

We consider a final type of integer sequence, the number of distinct values an invariant could obtain for a given order.
These sequences are not restricted to integer invariants. 
As an example, consider the integer sequence defined by the number of unique chromatic polynomials of graphs of a given order.

In this paper, we restrict the classes examined to those of simple connected graphs.
Unless otherwise stated, any referenced graph is assumed to be simple and connected.
Provided one had a means of enumeration, an extension of this program to other classes would be straightforward.
Exhaustive generation algorithms are known for many specialized classes such as bipartite graphs, digraphs, multigraphs, regular graphs, cubic graphs, snarks, trees and maximal triangle-free graphs\cite{mckay2014practical, meringer1999fast, brinkmann1996fast, brinkmann2013generation, sawada2006generating, brandt2000generation}. 

\section{Methods}

Using the \texttt{geng -c} command from \texttt{nauty} \cite{mckay2014practical}, we enumerated the the simple connected graphs up to order $n \le 10$.
Our calculations drew upon a large number of open source libraries and tools. 
Many of the graph invariant calculations were done with either \texttt{networkx} \cite{SciPyProceedings_11} or \texttt{graph-tool} \cite{Tiago2014graph}.
The invariants that were computable with integer or linear programming were done with PuLP \cite{mitchell2011pulp}.
For each graph we computed a series of invariants, which for completeness are described in Appendix \ref{app:invariants}.
A full table of all sequences submitted to the OEIS can be found in Appendix \ref{app:submittedseq}.
The relations between the invariants are numerous and are indexed online\cite{Travis2014Encyclopdia}.

Since the graphs we consider are loopless and undirected, the edge incidence information for each graph requires $n(n+1)/2$ bits of storage.
The largest we computed is of order 10, which requires 45 bits. 
This can be efficiently stored as a 64 bit unsigned integer using a binary representation.
Graphs of order 11 would also be possible to be stored in this representation, but graphs of order 12 would not (requiring 66 bits).
Internally, we used SQLite3 as the database back end for portability reasons.
Most of the invariants were stored as integers and compressed, the invariant database is about $850$ MB.
Some non-integer invariants, such as the Tutte polynomial, were stored in a separate, specialized database.

The creation of the database started with the full enumeration of the graphs themselves.
The calculation of the invariants was an embarrassingly parallel problem, as each graph was independent.
The invariants were stored in a denormalized database with each graph stored as a row and each column containing the value of a particular invariant for that graph. 
We placed an index on each of the invariant columns to allow for fast querying at the expense of disk space.

Sequences were constructed by first enumerating all the unique values for each invariant and applying a condition of inequality. 
For example, the sequence \OEIS{A241454} describes the graphs whose automorphism group has cardinality equal to two, while the sequence \OEIS{A086216} describes graphs with a vertex connectivity greater than or equal to three. 
Sequences can involve more than one invariant condition.
A sequence with only a single invariant condition is called a primary sequence, while a sequence involving two or more invariant conditions is called a secondary sequence.
For a sequence to be considered for submission to the OEIS, it must have at least four non-zero entries. 

The suggestive relations described in the previous section were built in a similar way.
For each pair of invariants conditions, the sets of graphs related to the sequences were compared for set equality, exclusivity, or as subsets.
In addition to providing useful information, some of these relations provided a check on the accuracy of the algorithms.
For example, one of the equality relations states that: a graph is a tree $\Leftrightarrow$ a graph has infinite girth.
This obviously follows from the fact that trees contain no cycles and the girth of a graph is the length of the shortest cycle (which is defined to be infinite for graphs without cycles).
Other examples are less trivial and hint at classic proofs.
For all simple connected graphs of order $n\le10$, it is true (and thus suggestive) that: a graph has a girth of five $\rightarrow$ a graph has chromatic number three. 
This eventually breaks down, in fact for a large enough $n$ there exists a graph for any given minimal girth and chromatic number\cite{erdos1959graph}.

The sequences counting distinct invariant values involved the creation of an additional database to store the non-integer invariants. 
Each of these special invariants required their own individualized table.
For example, the Tutte polynomial was stored as a denormalized table of graph identifiers with a row for each term of the polynomial containing the degree and coefficient.
Since this specialized database was much larger than the invariant integer database, it was not included in the public version.

\section{Discussion}

We have computed the Encyclopedia of Finite Graphs, an exhaustive enumeration of all simple connected graphs and their invariant values for graphs of order ten or less.
This database was utilized to investigate integer sequences among the invariants and supplement the OEIS.
While fully functional, there are several immediate extensions to the project that we aim to implement in the future.

From a technological standpoint, it is currently possible to extend to database to simple connected graphs of larger orders.
There are 1,006,700,565 and 164,059,830,476 graphs with 11 and 12 vertices respectively, a large but tractable size for local storage. 
Larger orders would require prohibitive storage requirements as the number of graphs grows quickly (sequence \OEIS{A001349}).
We choose a cutoff of order ten as a reasonable compromise between database size, distribution options and computational cost. 

There are many other classes of graphs to consider.
As mentioned earlier, not only are there other enumeration algorithms, but other projects have enumerated graphs of larger order with particular characteristics, such as bipartite, Eulerian, cubic, or certain regular graphs\cite{mckay2014combinatorial,royle2014small,meringer2014regular}.
An extension of the Encyclopedia could incorporate graphs that others have already determined along with a computation of their respective invariants.

In addition to expanding the volume of graphs considered, we can extend the number of computed invariants. 
Since the Encyclopedia code is openly available\cite{Travis2014Encyclopdia}, programs to compute additional invariants can be added by anyone. 
Large collections of graph invariants have been compiled (though not computed) under the Global Constraint Catalog, which also extensively details the relations among invariants\cite{beldiceanu2005global}.
As previously mentioned, the GraPHedron project aims to discover relationships among invariants using linear programming techniques. 
The Encyclopedia, with a greater number of graphs and invariants, could be employed to extend their results. 

\begin{figure}[p]
\begin{center}
  \includegraphics[width=0.75\textwidth]{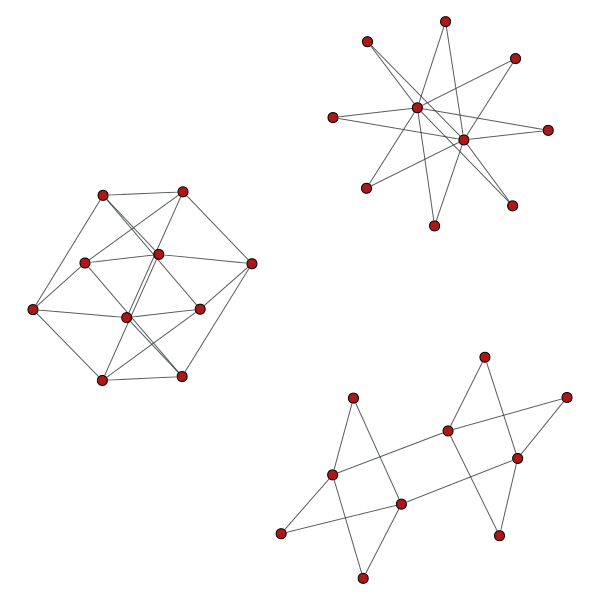}
\end{center}
  \caption{
    An example query to the Encyclopedia using specific invariant conditions. 
    The command \texttt{python viewer.py 10 -i is\_bipartite 1 -i is\_integral 1 -i is\_eulerian 1} displays the three graphs that are simultaneously bipartite, integral and Eulerian with ten vertices. 
  }
  \label{fig:examplequery}
\end{figure}

Being restricted to graphs of small order, the Encyclopedia's contribution is the compilation of graphs and their invariants into a queryable SQL database that is open-sourced and readily extensible. 
The Encyclopedia allows users to query specific invariant conditions and visualize the matching graphs.
As an example, there are exactly three graphs with 10 vertices that are simultaneously bipartite, integral and Eulerian. 
These graphs are illustrated in Figure \ref{fig:examplequery}.
As the Encyclopedia grows however, the published static database will be insufficient to handle these extensions.
We plan to create a web interface that would allow users to query graphs based on their name their names (e.g. ``Peterson graph") or their invariant conditions.

\clearpage 

\bibliographystyle{unsrt}
\bibliography{refs}

\newpage 

\begin{appendices}
\section{Invariant Descriptions}
\label{app:invariants}

The invariant descriptions are broadly organized into sections of algebraic graph theory, topological graph theory, and invariants related to structure-like connectivity.
Readers are instructed to look to standard references for more comprehensive definitions \cite{diestel2010graph, bondy2008graph, bollobas1998modern}.

\subsubsection*{Algebraic Graph Theory}

The \textit{adjacency matrix} $A$ has the value at $A_{ij}$ equal to the number of edges joining vertex $i$ to $j$. 
For simple graphs, this is a $\{0,1\}$ matrix with zeros down the diagonal. 
The \textit{characteristic polynomial} for a matrix is given by its eigenvalues, $\det(\lambda I - A)$.
A sequence of eigenvalues $\lambda_1, \lambda_2, \ldots$ is called the \textit{spectrum}.
A graph is $integral$ if all the values of the adjacency spectrum are integral, $\lambda_k \in \mathbb{Z}$. 
A graph is $real$ if $\lambda_k \in \mathbb{R}$.

An \textit{independent vertex/edge set} is a subset of vertices/edges of a graph such that no two vertices/edges in the subset share an edge/vertex. 
The \textit{maximal independent vertex set} or \textit{independence number}, is the cardinality of the largest independent vertex set.
A graph is \textit{bipartite} if there exists two disjoint independent subsets of the vertices whose union is $V$.
The \textit{maximal independent edge set} is the cardinality of the largest independent edge set.
The total number of independent edges sets is sometimes referred to as the number of \textit{edge matchings}, or the \textit{Hosoya index}\cite{hosoya1971topological}.
For graphs that are bounded in treewidth the Hosoya index is fixed-parameter tractable.
In general however, the Hosoya index $\mathcal{\#P}$-complete\cite{jerrum1987two}.

The \textit{automorphism group} is formed by all mappings of the vertices onto themselves which preserve isomorphism. 
The cardinality of this group is the \textit{automorphism number}.
Both \texttt{nauty} and \texttt{BLISS} \cite{junttila2007engineering,mckay2014practical} can compute the automorphism group and the number.

The \textit{Tutte polynomial} is a bivariate polynomial which encodes information about the graph's connectedness.
The polynomial $T_g(x,y)$ is defined via the recurrence relation, $T_g = T_{g-e} + T_{g/e}$ for any edge $e \in E$ that is not a loop or a bridge. 
Here $g-e$ denotes edge removal, $g / e$ denotes edge contraction and a bridge is an edge whose removal disconnects the graph.
A graph with only $i$ loops and $j$ bridges has the base case of $T_g(x,y) = x^i y^j$.
The \textit{chromatic polynomial} is a specialization of the Tutte polynomial 
\begin{equation}
\VARchromaticpolynomial{k} = (-1)^{|V|-c}k^{c}T_g(1-k,0)
\end{equation}
where $c$ is the number of connected components of the graph.
The chromatic polynomial $\VARchromaticpolynomial{k}$ gives the number of proper $k$-colorings on a graph. 
A proper $k$-coloring of a graph is an assignment of $k$ colors to each vertex so that no two adjacent vertices have the same color. 
The \textit{chromatic number} $\VARchromaticnumber$, is the smallest non-zero $k$ such that $\VARchromaticpolynomial{k}>0$.

The \textit{fractional coloring number} is motivated by the chromatic number.
A graph is \textit{$a:b$ colorable} if each vertex is assigned $b$ colors out of a palette of $a$ total colors such that the coloring of each adjacent vertex is disjoint. 
The smallest $a$ such that the graph can be assigned to $b$ colors is the $b$-fold chromatic number and is denoted by $\VARfractionalchromaticnumber{b}$. 
Note that the standard chromatic number is $\VARchromaticnumber{} = \VARfractionalchromaticnumber{1}$.
The fractional coloring number is defined to be
\begin{equation}
\VARfractionalchromaticnumber{f} = 
\lim_{b \rightarrow \infty} \frac{\VARfractionalchromaticnumber{b}}{b} = 
\inf _ {b} \frac{\VARfractionalchromaticnumber{b}}{b} 
\end{equation}
Not only does the limit exist, but $\VARfractionalchromaticnumber{f} \in \mathbb{Q}$\cite{scheinerman2011fractional}.
The integral chromatic number need not be equal to its fractional counter part, in fact 
$\VARfractionalchromaticnumber{f} \le \VARchromaticnumber$.
If the equality does not hold, we say the graph has a \textit{fractional chromatic gap}.
Additionally, we can formulate $\VARfractionalchromaticnumber{f}$ as a solution to a linear program.
Let $I(g)$ be the set of all independent vertex sets and $I(g,v)$ be the set of independent vertex sets that contain vertex $v$. 
Define a system of non-negative variables $v_i$.
Now $\VARfractionalchromaticnumber{f}$ is the minimum of $\sum_{I \in I(g)} v_i$ subject to the constraints $\sum_{I \in I(g,x)} v_i \ge 1$ on each $v_i$.

\subsubsection*{Topological Graph Theory}

The \textit{crossing number} is the minimum possible number of edge crossings for an embedding of the graph in a plane.
A graph is \textit{planar} if the crossing number is zero, in other words, a graph is planar if it admits an embedding in the plane with no overlapping edges.
A graph is \textit{toroidal} if the crossing number is one.
The \textit{genus} is the minimum number of edges that must be added to a graph to make it planar.
Due to a lack of freely available algorithms, planarity is the only embedding currently considered.

\subsubsection*{Connectivity, Cycles and Subgraphs}

The \textit{degree} of vertex $k$ is the number of edges incident to the vertex, $d(v_k) = \sum_j A_{k j}$. 
The \textit{degree sequence} is formed by listing the vertex degrees in descending order, and the number of distinct degree sequences for graphs for simple connected graphs is OEIS \OEIS{A007721}.
The degree matrix is a diagonal matrix defined by the degree sequence.
The \textit{Laplacian matrix} is the difference of the diagonal matrix and the adjacency matrix.
The \textit{Laplacian polynomial} is the characteristic polynomial of the Laplacian matrix.
For connected graphs the spectrum of the Laplacian matrix is real and $\lambda_i \le 0, \lambda_0 = 0$.
The second-largest eigenvalue is called the \textit{spectral gap}.

Let $D$ denote the graph distance matrix, where $D_{ij}$ is the shortest path (geodesic) from vertex $i$ to $j$.
A graph is \textit{connected} if there is a path between all pairs of vertices.
If there is no path between vertices $i$ and $j$, $D_{ij}=\infty$.
By construction, all distances in connected graphs are finite.
The \textit{eccentricity} of a vertex $k$ is the maximum value of $D_{k j}$ for all $j$.
The \textit{radius}/\textit{diameter} is the minimum/maximum value of eccentricity for all vertices.

A \textit{$k$-regular graph} is one in which every vertex has degree $k$. 
A regular graph is \textit{strongly regular}, if there exist integers $a$ and $b$ such that any two adjacent vertices have $a$ neighbors in common and any two non-adjacent vertices have $b$ neighbors in common. 
A graph is \textit{distance regular} if two conditions are met. 
First, for all pairs of vertices $v, w$ and any pair of integers $i,j = 0,1,\ldots,d$ where $d$ is the graph diameter, the number of vertices from $v$ at a distance $i$ and the number of vertices from $w$ at a distance $j$ depend only on $i, j$.
Secondly, if the distance between $v$ and $w$ is independent of the choice of $v$ and $w$\cite{brouwer2012distance}.

A \textit{circuit} is a path that begins and ends at the same vertex.
A graph is \textit{Hamiltonian}/\textit{Eulerian} if there is a circuit that visits each vertex/edge exactly once.

A vertex is an \textit{articulation point} if its removal disconnects the graph. 
The \textit{vertex/edge connectivity} of a graph is the minimum number of vertices/edges needed to disconnect the graph.
A vertex joined to only a single edge is an \textit{end point}.

The \textit{cycle space} of a graph is the set of all of its Eulerian subgraphs, where every member can be constructed as the symmetric difference of members of the \textit{cycle basis}. 
A \textit{tree} is an acyclic connected graph whose cycle basis is the empty set.
The length of the shortest/longest cycle in a graph is its \textit{girth}/\textit{circumference}.
Trees are defined to have infinite girth and circumference.
Due to technical limitations we record this value as a zero which is unambiguous for connected loopless graphs.
A \textit{chord} of a cycle is an edge with one vertex belonging to the cycle and one edge outside of the cycle.
In a \textit{chordal graph}, all cycles of order four or more have a cycle chord.   
 
A \textit{vertex subgraph} or simply a subgraph, is a subset of vertices and the edges which are common to all members of this subset. 
Let $\VARsubgraph(g,h)$ equal the number of vertex sets that form subgraphs of $h$ in $g$.
The \textit{maximum clique number} is the largest non-zero value of $\VARsubgraph(g, K_n)$, where $K_n$ is the complete graph on $n$ vertices. 
A graph is \textit{triangle free} or \textit{square free} if $\VARsubgraph(g,C_3)=0$ or $\VARsubgraph(g,C_4)=0$ respectively, where $C_n$ is the cycle graph on $n$ vertices.
In addition, we have checked for the subgraphs $K_4, K_5, C_5, C_6$ and several named graphs.
Illustrated below from left-to-right are the 
the bull graph $\subgraphBULL{}$,
the bowtie graph $\subgraphBOWTIE{}$, 
the open-bowtie graph $\subgraphOPENBOWTIE{}$,
and the diamond graph $\subgraphDIAMOND{}$.
\begin{center}
  \includegraphics[width=0.5\textwidth]{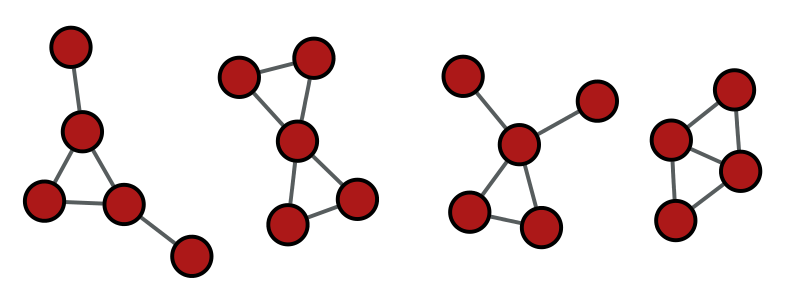}
\end{center}

\clearpage 

\section{Submitted sequences}
\label{app:submittedseq}

The sequences that were already present in the OEIS but extended by this project are listed in Section \ref{sec:seq_extended}.
The distinct sequences were those formed by considering the number of unique values possible at a fixed order and are listed in Section \ref{sec:seq_distinct}.
Sequences involving a single invariant, the primary sequences, are listed in Section \ref{sec:seq_primary}.
The secondary sequences involve two invariants and are listed in Section \ref{sec:seq_secondary}.

Each sequence was hand-checked against the OEIS.
We validated that either the sequence was unique or the invariant description matched the sequence found.
We reiterate that all results shown are for simple, connected graphs. 
For some invariants the definition may be ambiguous for small orders.
When possible, we tried to refer to conventions found in the literature and the OEIS for these cases.
For brevity the symbols describing the invariants are shown below.

\begin{longtable}{@{\extracolsep{\fill}} l l}
\toprule
Invariant symbol & Description \\
\midrule
$\VARautomorphismgroupn$ & Automorphism group size \\
$\VARchromaticnumber$ & Chromatic number \\
$\VARvertexconnectivity$, $\VARedgeconnectivity$ & vertex/edge connectivity \\
$\VARdiameter, \VARgirth$ & diameter, girth \\
$\VARnarticulationpoints$ & number of articulation points \\
$\VARmaximalindependentvertexset, \VARmaximalindependentedgeset$ & size of maximal independent vertex/edge set \\
$\namedsubgraph{h}$ & number of subgraphs of $h$ in $g$ \\
$\VARtree,\VARhamiltonian,\VAReulerian, \VARintegral, \VARrealspectrum$ & 
properties: tree, Hamiltonian, Eulerian, integral, real \\
$\VARkregular,\VARstronglyregular,\VARdistanceregular, \VARbipartite$ & 
properties: regular, strongly regular, distance regular, bipartite \\
$\VARchordal, \VARfractionaldualitygapvertexchromatic$  &
properties: chordal, chromatic gap
\end{longtable}

\newenvironment{invariantTable}[2]
{
  \begin{longtable}{@{\extracolsep{\fill}} #1}
    \toprule
    OEIS & Invariant & #2 \\
    \midrule\bottomrule
}
{
  \end{longtable}
}

\subsection{Extended sequences}
\label{sec:seq_extended}

\begin{invariantTable}{l l r r r r r r r r r r}
{\multicolumn{10}{l}{Sequence: $\SEQ_1, \SEQ_2, \ldots, \SEQ_{10}$, terms in $(\cdot)$ were added by this project}}
\OEIS{A086216} & $\VARvertexconnectivity > 3$ & 0 & 0 & 0 & 0 & 1 & 4 & 25 & 384 & 14480 & (1211735) \\
\OEIS{A086217} & $\VARvertexconnectivity > 4$ & 0 & 0 & 0 & 0 & 0 & 1 & 4 & 39 & 1051 & (102630) \\
\OEIS{A052446} & $\VARedgeconnectivity =2$ & 	0& 1& 1& 3& 10& 52& 351& 3714& 63638& (1912203) \\
\OEIS{A052447} & $\VARedgeconnectivity =3$ & 0& 0& 1& 2& 8& 41& 352& (4820) & (113256) & (4602039) \\
\OEIS{A052448} & $\VARedgeconnectivity =4$ & 	0& 0& 0& 1& 2& 15& 121& (2159) & (68715) & (3952378) \\
\OEIS{A088741} & $\VARstronglyregular(g)=1$ & 	1& 1& 1& 2& 2& 3& 1& 3& 3& (5) \\
\end{invariantTable}

\clearpage 

\subsection{Distinct sequences}
\label{sec:seq_distinct}

\begin{invariantTable}{l l r r r r r r r r r r}
{\multicolumn{10}{l}{Sequence: $\SEQ_1, \SEQ_2, \ldots, \SEQ_{10}$}}
\OEIS{A245881} & Tutte poly & 1 & 1 & 2 & 5 & 16 & 73 & 532 & 7245 & 192605 & 9717156 \\
\OEIS{A245883} & chromatic poly & 1 & 1 & 2 & 5 & 14 & 50 & 231 & 1650 & 21121 & 584432 \\
\OEIS{A245880} & characteristic poly & 1 & 1 & 2 & 6 & 21 & 111 & 821 & 10423 & 236064 & 10375796 \\
\OEIS{A245882} & Laplacian poly & 1 & 1 & 2 & 6 & 21 & 110 & 793 & 10251 & 239307 & 10985229 \\
\OEIS{A245879} & fractional chromatic& 1 & 1 & 2 & 3 & 5 & 7 & 11 & 17 & 29 & 50 \\
\end{invariantTable}

\subsection{Primary sequences}
\label{sec:seq_primary}

\begin{invariantTable}{l l r r r r r r r r r r}
{\multicolumn{10}{l}{Sequence: $\SEQ_1, \SEQ_2, \ldots, \SEQ_{10}$}}
\OEIS{A241454} & $\VARautomorphismgroupn =2$ & 0 & 1 & 1 & 2 & 9 & 37 & 317 & 4098 & 84602 & 2933996 \\
\OEIS{A241455} & $\VARautomorphismgroupn =4$ & 0 & 0 & 0 & 1 & 3 & 28 & 198 & 1971 & 29047 & 672516 \\
\OEIS{A241456} & $\VARautomorphismgroupn =6$ & 0 & 0 & 1 & 1 & 1 & 7 & 31 & 221 & 3025 & 68033 \\
\OEIS{A241457} & $\VARautomorphismgroupn =8$ & 0 & 0 & 0 & 1 & 2 & 9 & 55 & 499 & 6017 & 107312 \\
\OEIS{A241458} & $\VARautomorphismgroupn =10$ & 0 & 0 & 0 & 0 & 1 & 1 & 1 & 3 & 13 & 123 \\
\OEIS{A241459} & $\VARautomorphismgroupn =12$ & 0 & 0 & 0 & 0 & 3 & 10 & 51 & 356 & 3395 & 49862 \\
\OEIS{A241460} & $\VARautomorphismgroupn =14$ & 0 & 0 & 0 & 0 & 0 & 0 & 2 & 2 & 2 & 6 \\
\OEIS{A241461} & $\VARautomorphismgroupn =16$ & 0 & 0 & 0 & 0 & 0 & 3 & 10 & 123 & 992 & 14026 \\
\OEIS{A241462} & $\VARautomorphismgroupn =20$ & 0 & 0 & 0 & 0 & 0 & 0 & 2 & 6 & 29 & 199 \\
\OEIS{A241463} & $\VARautomorphismgroupn =24$ & 0 & 0 & 0 & 1 & 1 & 1 & 14 & 118 & 1247 & 17191 \\
\OEIS{A241464} & $\VARautomorphismgroupn =36$ & 0 & 0 & 0 & 0 & 0 & 1 & 3 & 16 & 132 & 1341 \\
\OEIS{A241465} & $\VARautomorphismgroupn =48$ & 0 & 0 & 0 & 0 & 0 & 4 & 14 & 65 & 504 & 5215 \\
\OEIS{A241466} & $\VARautomorphismgroupn =72$ & 0 & 0 & 0 & 0 & 0 & 1 & 2 & 16 & 124 & 1070 \\
\OEIS{A241467} & $\VARautomorphismgroupn =120$ & 0 & 0 & 0 & 0 & 1 & 1 & 1 & 5 & 21 & 211 \\
\OEIS{A241468} & $\VARautomorphismgroupn =144$ & 0 & 0 & 0 & 0 & 0 & 0 & 3 & 12 & 51 & 477 \\
\OEIS{A241469} & $\VARautomorphismgroupn =240$ & 0 & 0 & 0 & 0 & 0 & 0 & 3 & 8 & 51 & 336 \\
\OEIS{A241470} & $\VARautomorphismgroupn =720$ & 0 & 0 & 0 & 0 & 0 & 1 & 1 & 4 & 13 & 60 \\
\OEIS{A241471} & $\VARautomorphismgroupn =5040$ & 0 & 0 & 0 & 0 & 0 & 0 & 1 & 1 & 1 & 5 \\
\OEIS{A241702} & $\VARchromaticnumber =7$ & 0 & 0 & 0 & 0 & 0 & 0 & 1 & 6 & 110 & 4125 \\
\OEIS{A241703} & $\VARedgeconnectivity =4$ & 0 & 0 & 0 & 0 & 1 & 3 & 25 & 378 & 14306 & 1141575 \\
\OEIS{A241704} & $\VARedgeconnectivity =5$ & 0 & 0 & 0 & 0 & 0 & 1 & 3 & 41 & 1095 & 104829 \\
\OEIS{A241705} & $\VARedgeconnectivity =6$ & 0 & 0 & 0 & 0 & 0 & 0 & 1 & 4 & 65 & 3441 \\
\OEIS{A241706} & $\VARdiameter =2$ & 0 & 0 & 1 & 4 & 14 & 59 & 373 & 4154 & 91518 & 4116896 \\
\OEIS{A241707} & $\VARdiameter =3$ & 0 & 0 & 0 & 1 & 5 & 43 & 387 & 5797 & 148229 & 6959721 \\
\OEIS{A241708} & $\VARdiameter =4$ & 0 & 0 & 0 & 0 & 1 & 8 & 82 & 1027 & 19320 & 598913 \\
\OEIS{A241709} & $\VARdiameter =5$ & 0 & 0 & 0 & 0 & 0 & 1 & 9 & 125 & 1818 & 37856 \\
\OEIS{A241710} & $\VARdiameter =6$ & 0 & 0 & 0 & 0 & 0 & 0 & 1 & 12 & 180 & 2928 \\
\OEIS{A241711} & $\VARgirth =3$ & 0 & 0 & 1 & 3 & 15 & 93 & 792 & 10833 & 259420 & 11704309 \\
\OEIS{A241712} & $\VARgirth =4$ & 0 & 0 & 0 & 1 & 2 & 11 & 43 & 234 & 1498 & 11451 \\
\OEIS{A241713} & $\VARgirth =5$ & 0 & 0 & 0 & 0 & 1 & 1 & 5 & 18 & 82 & 539 \\
\OEIS{A241714} & $\VARgirth =6$ & 0 & 0 & 0 & 0 & 0 & 1 & 1 & 7 & 25 & 137 \\
\OEIS{A241715} & $\VARgirth =7$ & 0 & 0 & 0 & 0 & 0 & 0 & 1 & 1 & 6 & 20 \\
\OEIS{A241767} & $\VARnarticulationpoints =1$ & 0 & 0 & 1 & 2 & 7 & 33 & 244 & 2792 & 52448 & 1690206 \\
\OEIS{A241768} & $\VARnarticulationpoints =2$ & 0 & 0 & 0 & 1 & 3 & 17 & 101 & 890 & 11468 & 239728 \\
\OEIS{A241769} & $\VARnarticulationpoints =3$ & 0 & 0 & 0 & 0 & 1 & 5 & 32 & 242 & 2461 & 35839 \\
\OEIS{A241770} & $\VARnarticulationpoints =4$ & 0 & 0 & 0 & 0 & 0 & 1 & 7 & 60 & 527 & 6056 \\
\OEIS{A241771} & $\VARnarticulationpoints =5$ & 0 & 0 & 0 & 0 & 0 & 0 & 1 & 9 & 97 & 1029 \\
\OEIS{A241782} & $\VARissubgraphfreeKfive =0$ & 1 & 1 & 2 & 6 & 20 & 107 & 802 & 10252 & 232850 & 9905775 \\
\OEIS{A241784} & $\VARissubgraphfreeCfive =0$ & 1 & 1 & 2 & 6 & 13 & 44 & 144 & 577 & 2457 & 12499 \\
\OEIS{A242790} & $\VARissubgraphfreediamond =0$ & 1 & 1 & 2 & 4 & 11 & 39 & 165 & 967 & 7684 & 87012 \\
\OEIS{A242792} & $\VARissubgraphfreebowtie =0$ & 1 & 1 & 2 & 6 & 15 & 60 & 273 & 1769 & 14836 & 174111 \\
\OEIS{A242791} & $\VARissubgraphfreeopenbowtie =0$ & 1 & 1 & 2 & 6 & 11 & 34 & 98 & 408 & 1957 & 12740 \\
\OEIS{A243243} & $\VARissubgraphfreeCfour >0$ & 0 & 0 & 0 & 3 & 13 & 93 & 796 & 10931 & 260340 & 11713182 \\
\OEIS{A243246} & $\VARissubgraphfreeCfive >0$ & 0 & 0 & 0 & 0 & 8 & 68 & 709 & 10540 & 258623 & 11704072 \\
\OEIS{A243245} & $\VARissubgraphfreeKthree >0$ & 0 & 0 & 1 & 3 & 15 & 93 & 794 & 10850 & 259700 & 11706739 \\
\OEIS{A243244} & $\VARissubgraphfreeKfour >0$ & 0 & 0 & 0 & 1 & 4 & 30 & 317 & 5511 & 165165 & 8932499 \\
\OEIS{A243242} & $\VARissubgraphfreeKfive >0$ & 0 & 0 & 0 & 0 & 1 & 5 & 51 & 865 & 28230 & 1810796 \\
\OEIS{A243250} & $\VARissubgraphfreediamond >0$ & 0 & 0 & 0 & 2 & 10 & 73 & 688 & 10150 & 253396 & 11629559 \\
\OEIS{A243248} & $\VARissubgraphfreebull >0$ & 0 & 0 & 0 & 0 & 12 & 86 & 773 & 10777 & 259390 & 11705139 \\
\OEIS{A243249} & $\VARissubgraphfreebowtie >0$ & 0 & 0 & 0 & 0 & 6 & 52 & 580 & 9348 & 246244 & 11542460 \\
\OEIS{A243247} & $\VARissubgraphfreeopenbowtie >0$ & 0 & 0 & 0 & 0 & 10 & 78 & 755 & 10709 & 259123 & 11703831 \\
\OEIS{A241814} & $\VARisdistanceregular =1$ & 1 & 1 & 1 & 2 & 2 & 4 & 2 & 5 & 4 & 7 \\
\OEIS{A241839} & $\VARiskregular =0$ & 1 & 0 & 1 & 4 & 19 & 107 & 849 & 11100 & 261058 & 11716404 \\
\OEIS{A241840} & $\VARisdistanceregular =0$ & 0 & 0 & 1 & 4 & 19 & 108 & 851 & 11112 & 261076 & 11716564 \\
\OEIS{A241841} & $\VARistree =0$ & 0 & 0 & 1 & 4 & 18 & 106 & 842 & 11094 & 261033 & 11716465 \\
\OEIS{A241842} & $\VARisintegral =0$ & 0 & 0 & 1 & 4 & 18 & 106 & 846 & 11095 & 261056 & 11716488 \\
\OEIS{A241843} & $\VARischordal =0$ & 0 & 0 & 0 & 1 & 6 & 54 & 581 & 9503 & 249169 & 11607032 \\
\OEIS{A242952} & $\VARisrealspectrum =1$ & 1 & 1 & 1 & 3 & 11 & 54 & 539 & 7319 & 209471 & 10000304 \\
\OEIS{A242953} & $\VARisrealspectrum =0$ & 0 & 0 & 1 & 3 & 10 & 58 & 314 & 3798 & 51609 & 1716267 \\
\OEIS{A243241} & $\VARisstronglyregular =0$ & 0 & 0 & 1 & 4 & 19 & 109 & 852 & 11114 & 261077 & 11716566 \\
\OEIS{A243251} & $\VARhasfractionaldualitygapvertexchromatic =1$ & 0 & 0 & 0 & 0 & 1 & 3 & 33 & 496 & 16464 & 969293 \\
\OEIS{A243252} & $\VARhasfractionaldualitygapvertexchromatic =0$ & 1 & 1 & 2 & 6 & 20 & 109 & 820 & 10621 & 244616 & 10747278 \\
\OEIS{A243781} & $\VARmaximalindependentvertexset =2$ & 0 & 0 & 1 & 4 & 11 & 34 & 103 & 405 & 1892 & 12166 \\
\OEIS{A243782} & $\VARmaximalindependentvertexset =3$ & 0 & 0 & 0 & 1 & 8 & 63 & 524 & 5863 & 100702 & 2880002 \\
\OEIS{A243783} & $\VARmaximalindependentvertexset =4$ & 0 & 0 & 0 & 0 & 1 & 13 & 205 & 4308 & 135563 & 7161399 \\
\OEIS{A243784} & $\VARmaximalindependentvertexset =5$ & 0 & 0 & 0 & 0 & 0 & 1 & 19 & 513 & 21782 & 1576634 \\
\OEIS{A243800} & $\VARmaximalindependentedgeset =2$ & 0 & 0 & 0 & 5 & 20 & 16 & 22 & 29 & 37 & 46 \\
\OEIS{A243801} & $\VARmaximalindependentedgeset =3$ & 0 & 0 & 0 & 0 & 0 & 95 & 830 & 790 & 1479 & 2625 \\
\end{invariantTable}

\clearpage 

\subsection{Secondary sequences}
\label{sec:seq_secondary}

\begin{invariantTable}{l l l r r r r r r r r r r}
{& \multicolumn{10}{l}{Sequence: $\SEQ_1, \SEQ_2, \ldots, \SEQ_{10}$}}
\OEIS{A243270} & $\VARishamiltonian =1$, & $\VARisbipartite =1$ & 1 & 0 & 0 & 1 & 0 & 4 & 0 & 24 & 0 & 473 \\
\OEIS{A243271} & $\VARishamiltonian =1$, & $\VARisdistanceregular =1$ & 1 & 0 & 1 & 2 & 2 & 4 & 2 & 5 & 4 & 6 \\
\OEIS{A243272} & $\VARishamiltonian =1$, & $\VARiseulerian =1$ & 1 & 0 & 1 & 1 & 2 & 5 & 21 & 120 & 1312 & 26525 \\
\OEIS{A243273} & $\VARishamiltonian =1$, & $\VARisintegral =0$ & 0 & 0 & 0 & 1 & 7 & 43 & 379 & 6185 & 177071 & 9305068 \\
\OEIS{A243274} & $\VARishamiltonian =1$, & $\VARisintegral =1$ & 1 & 0 & 1 & 2 & 1 & 5 & 4 & 11 & 12 & 50 \\
\OEIS{A243275} & $\VARishamiltonian =1$, & $\VARissubgraphfreeKthree =0$ & 1 & 0 & 0 & 1 & 1 & 4 & 5 & 35 & 130 & 1293 \\
\OEIS{A243276} & $\VARishamiltonian =1$, & $\VARissubgraphfreeKfour =0$ & 1 & 0 & 1 & 2 & 5 & 29 & 188 & 2481 & 52499 & 1857651 \\
\OEIS{A243319} & $\VARisbipartite =1$, & $\VARisdistanceregular =1$ & 1 & 1 & 0 & 1 & 0 & 2 & 0 & 3 & 0 & 3 \\
\OEIS{A243320} & $\VARisbipartite =1$, & $\VARiseulerian =1$ & 1 & 0 & 0 & 1 & 0 & 2 & 1 & 6 & 7 & 29 \\
\OEIS{A243321} & $\VARisbipartite =1$, & $\VARisplanar =1$ & 1 & 1 & 1 & 3 & 5 & 16 & 41 & 158 & 582 & 2749 \\
\OEIS{A243322} & $\VARisdistanceregular =1$, & $\VARiseulerian =1$ & 1 & 0 & 1 & 1 & 2 & 2 & 2 & 3 & 4 & 4 \\
\OEIS{A243323} & $\VARisintegral =0$, & $\VARisbipartite =1$ & 0 & 0 & 1 & 2 & 4 & 14 & 43 & 179 & 730 & 4019 \\
\OEIS{A243324} & $\VARisintegral =0$, & $\VARiseulerian =1$ & 0 & 0 & 0 & 0 & 2 & 6 & 33 & 180 & 1773 & 31006 \\
\OEIS{A243325} & $\VARisintegral =0$, & $\VARisplanar =1$ & 0 & 0 & 1 & 4 & 18 & 95 & 642 & 5962 & 71876 & 1052786 \\
\OEIS{A243326} & $\VARisintegral =0$, & $\VARissubgraphfreeKthree =0$ & 0 & 0 & 1 & 2 & 5 & 16 & 58 & 264 & 1380 & 9818 \\
\OEIS{A243327} & $\VARisintegral =0$, & $\VARissubgraphfreeKfour =0$ & 0 & 0 & 1 & 4 & 15 & 77 & 531 & 5597 & 95900 & 2784034 \\
\OEIS{A243328} & $\VARisintegral =1$, & $\VARisbipartite =1$ & 1 & 1 & 0 & 1 & 1 & 3 & 1 & 3 & 0 & 13 \\
\OEIS{A243329} & $\VARisintegral =1$, & $\VARisdistanceregular =1$ & 1 & 1 & 1 & 2 & 1 & 4 & 1 & 4 & 3 & 6 \\
\OEIS{A243330} & $\VARisintegral =1$, & $\VARiseulerian =1$ & 1 & 0 & 1 & 1 & 2 & 2 & 4 & 4 & 9 & 20 \\
\OEIS{A243331} & $\VARisintegral =1$, & $\VARisplanar =1$ & 1 & 1 & 1 & 2 & 2 & 4 & 4 & 12 & 9 & 19 \\
\OEIS{A243332} & $\VARisintegral =1$, & $\VARissubgraphfreeKthree =0$ & 1 & 1 & 0 & 1 & 1 & 3 & 1 & 3 & 0 & 14 \\
\OEIS{A243333} & $\VARisintegral =1$, & $\VARissubgraphfreeKfour =0$ & 1 & 1 & 1 & 1 & 2 & 5 & 5 & 9 & 15 & 38 \\
\OEIS{A243334} & $\VARissubgraphfreeKthree =0$, & $\VARisdistanceregular =1$ & 1 & 1 & 0 & 1 & 1 & 2 & 1 & 3 & 1 & 4 \\
\OEIS{A243335} & $\VARissubgraphfreeKthree =0$, & $\VARiseulerian =1$ & 1 & 0 & 0 & 1 & 1 & 2 & 3 & 8 & 19 & 62 \\
\OEIS{A243336} & $\VARissubgraphfreeKfour =0$, & $\VARiseulerian =1$ & 1 & 0 & 1 & 1 & 3 & 6 & 22 & 93 & 656 & 7484 \\
\OEIS{A243337} & $\VARissubgraphfreeKfour =0$, & $\VARisplanar =1$ & 1 & 1 & 2 & 5 & 17 & 79 & 478 & 4123 & 46666 & 648758 \\
\OEIS{A243338} & $\VARistree =1$, & $\VARisintegral =0$ & 0 & 0 & 1 & 2 & 2 & 5 & 10 & 23 & 47 & 105 \\
\OEIS{A243339} & $\VARissubgraphfreeKfour =0$, & $\VARisdistanceregular =1$ & 1 & 1 & 1 & 1 & 1 & 3 & 1 & 3 & 3 & 4 \\
\OEIS{A243545} & $\VARishamiltonian =1$, & $\VARissubgraphfreebowtie =0$ & 1 & 0 & 1 & 3 & 3 & 14 & 50 & 390 & 3627 & 52858 \\
\OEIS{A243546} & $\VARissubgraphfreebowtie =0$, & $\VARisdistanceregular =1$ & 1 & 1 & 1 & 2 & 1 & 2 & 1 & 3 & 1 & 4 \\
\OEIS{A243547} & $\VARissubgraphfreebowtie =0$, & $\VARiseulerian =1$ & 1 & 0 & 1 & 1 & 2 & 4 & 8 & 35 & 115 & 629 \\
\OEIS{A243548} & $\VARissubgraphfreebowtie =0$, & $\VARisintegral =1$ & 1 & 1 & 1 & 2 & 2 & 4 & 1 & 8 & 1 & 19 \\
\OEIS{A243549} & $\VARissubgraphfreebowtie =0$, & $\VARisintegral =0$ & 0 & 0 & 1 & 4 & 13 & 56 & 272 & 1761 & 14835 & 174092 \\
\OEIS{A243550} & $\VARissubgraphfreebowtie =0$, & $\VARisplanar =1$ & 1 & 1 & 2 & 6 & 15 & 58 & 255 & 1510 & 10766 & 94109 \\
\OEIS{A243551} & $\VARissubgraphfreebowtie =0$, & $\VARissubgraphfreeKfour =0$ & 1 & 1 & 2 & 5 & 14 & 56 & 256 & 1656 & 13952 & 163878 \\
\OEIS{A243552} & $\VARissubgraphfreebowtie =0$, & $\VARissubgraphfreebull =0$ & 1 & 1 & 2 & 6 & 8 & 25 & 77 & 333 & 1668 & 11355 \\
\OEIS{A243553} & $\VARishamiltonian =1$, & $\VARissubgraphfreebull =0$ & 1 & 0 & 1 & 3 & 1 & 4 & 5 & 35 & 130 & 1293 \\
\OEIS{A243554} & $\VARissubgraphfreebull =0$, & $\VARisdistanceregular =1$ & 1 & 1 & 1 & 2 & 1 & 2 & 1 & 3 & 1 & 4 \\
\OEIS{A243555} & $\VARissubgraphfreebull =0$, & $\VARiseulerian =1$ & 1 & 0 & 1 & 1 & 2 & 3 & 5 & 14 & 30 & 95 \\
\OEIS{A243556} & $\VARissubgraphfreebull =0$, & $\VARisintegral =1$ & 1 & 1 & 1 & 2 & 1 & 3 & 2 & 3 & 0 & 14 \\
\OEIS{A243557} & $\VARissubgraphfreebull =0$, & $\VARisintegral =0$ & 0 & 0 & 1 & 4 & 8 & 23 & 78 & 337 & 1690 & 11418 \\
\OEIS{A243558} & $\VARissubgraphfreebull =0$, & $\VARisplanar =1$ & 1 & 1 & 2 & 6 & 9 & 25 & 76 & 302 & 1360 & 7606 \\
\OEIS{A243559} & $\VARissubgraphfreebull =0$, & $\VARissubgraphfreeKfour =0$ & 1 & 1 & 2 & 5 & 9 & 26 & 80 & 340 & 1690 & 11432 \\
\OEIS{A243560} & $\VARishamiltonian =1$, & $\VARissubgraphfreediamond =0$ & 1 & 0 & 1 & 1 & 2 & 9 & 27 & 190 & 1750 & 25658 \\
\OEIS{A243561} & $\VARissubgraphfreediamond =0$, & $\VARisdistanceregular =1$ & 1 & 1 & 1 & 1 & 1 & 2 & 1 & 3 & 2 & 4 \\
\OEIS{A243562} & $\VARissubgraphfreediamond =0$, & $\VARiseulerian =1$ & 1 & 0 & 1 & 1 & 2 & 3 & 8 & 21 & 79 & 334 \\
\OEIS{A243563} & $\VARissubgraphfreediamond =0$, & $\VARisintegral =0$ & 0 & 0 & 1 & 3 & 10 & 35 & 162 & 964 & 7682 & 86994 \\
\OEIS{A243564} & $\VARissubgraphfreediamond =0$, & $\VARisintegral =1$ & 1 & 1 & 1 & 1 & 1 & 4 & 3 & 3 & 2 & 18 \\
\OEIS{A243565} & $\VARissubgraphfreediamond =0$, & $\VARisplanar =1$ & 1 & 1 & 2 & 4 & 11 & 38 & 159 & 882 & 6242 & 55316 \\
\OEIS{A243566} & $\VARissubgraphfreediamond =0$, & $\VARissubgraphfreeKfour =0$ & 1 & 1 & 2 & 4 & 11 & 39 & 165 & 967 & 7684 & 87012 \\
\OEIS{A243567} & $\VARissubgraphfreediamond =0$, & $\VARissubgraphfreebowtie =0$ & 1 & 1 & 2 & 4 & 10 & 36 & 141 & 784 & 5626 & 56249 \\
\OEIS{A243568} & $\VARissubgraphfreediamond =0$, & $\VARissubgraphfreebull =0$ & 1 & 1 & 2 & 4 & 9 & 26 & 80 & 340 & 1690 & 11432 \\
\OEIS{A243789} & $\VARissubgraphfreeopenbowtie =0$, & $\VARissubgraphfreediamond =0$ & 1 & 1 & 2 & 4 & 9 & 30 & 89 & 379 & 1864 & 12365 \\
\OEIS{A243790} & $\VARissubgraphfreeopenbowtie =0$, & $\VARishamiltonian =1$ & 1 & 0 & 1 & 3 & 3 & 9 & 13 & 59 & 203 & 1651 \\
\OEIS{A243791} & $\VARissubgraphfreeopenbowtie =0$, & $\VARiseulerian =1$ & 1 & 0 & 1 & 1 & 1 & 2 & 3 & 8 & 19 & 62 \\
\OEIS{A243792} & $\VARissubgraphfreeopenbowtie =0$, & $\VARisintegral =1$ & 1 & 1 & 1 & 2 & 1 & 4 & 1 & 3 & 0 & 15 \\
\OEIS{A243783} & $\VARissubgraphfreeopenbowtie =0$, & $\VARisintegral =0$ & 0 & 0 & 1 & 4 & 10 & 30 & 97 & 405 & 1957 & 12725 \\
\OEIS{A243794} & $\VARissubgraphfreeopenbowtie =0$, & $\VARisplanar =1$ & 1 & 1 & 2 & 6 & 11 & 33 & 94 & 370 & 1627 & 8895 \\
\OEIS{A243795} & $\VARissubgraphfreeopenbowtie =0$, & $\VARissubgraphfreebull =0$ & 1 & 1 & 2 & 6 & 7 & 22 & 65 & 285 & 1442 & 10106 \\
\OEIS{A243253} & $\VARischordal =1$, & $\VARiseulerian =1$ & 1 & 0 & 1 & 0 & 3 & 2 & 13 & 18 & 116 & 366 \\
\OEIS{A243785} & $\VARischordal =1$, & $\VARisintegral =0$ & 0 & 0 & 1 & 4 & 12 & 56 & 267 & 1605 & 11909 & 109525 \\
\OEIS{A243786} & $\VARischordal =1$, & $\VARisintegral =1$ & 1 & 1 & 1 & 1 & 3 & 2 & 5 & 9 & 2 & 14 \\
\OEIS{A243787} & $\VARischordal =1$, & $\VARisplanar =1$ & 1 & 1 & 2 & 5 & 14 & 52 & 228 & 1209 & 7463 & 52520 \\
\OEIS{A243788} & $\VARischordal =1$, & $\VARissubgraphfreeKfour =0$ & 1 & 1 & 2 & 4 & 11 & 35 & 124 & 500 & 2224 & 10640 \\
\OEIS{A243796} & $\VARishamiltonian =1$, & $\VARischordal =1$ & 1 & 0 & 1 & 2 & 4 & 15 & 58 & 360 & 2793 & 28761 \\
\OEIS{A243797} & $\VARissubgraphfreebowtie =0$, & $\VARischordal =1$ & 1 & 1 & 2 & 5 & 10 & 27 & 70 & 206 & 613 & 1942 \\
\OEIS{A243798} & $\VARissubgraphfreebull =0$, & $\VARischordal =1$ & 1 & 1 & 2 & 5 & 6 & 12 & 25 & 55 & 126 & 304 \\
\OEIS{A243799} & $\VARissubgraphfreeopenbowtie =0$, & $\VARischordal =1$ & 1 & 1 & 2 & 5 & 6 & 13 & 25 & 58 & 130 & 316 \\
\end{invariantTable}

\end{appendices}

\end{document}